\documentclass[12pt]{amsart}
\usepackage{amssymb,amsmath,amsthm}
\usepackage{amscd}
\input amssym.def
\input amssym
\newtheorem{theorem}{Theorem}[section]

\newtheorem{lemma}{Lemma}[section]

\newtheorem{remark}{Remark}[section]
\newtheorem{defi}{Definition}[section]
\newtheorem{prop}{Proposition}[section]

\newcommand{\be}{\begin{equation}}
\newcommand{\ee}{\end{equation}}

\renewcommand{\theequation}{\thesection.\arabic{equation}}
\renewcommand{\thetheorem}{\thesection.\arabic{theorem}}
\numberwithin{equation}{section}

\begin{document}

\title[] {A recursion identity for formal iterated logarithms and
  iterated exponentials}

\author{Thomas J. Robinson}


\begin{abstract}
We prove a recursive identity involving formal iterated logarithms and
formal iterated exponentials.  These iterated logarithms and
exponentials appear in a natural extension of the logarithmic formal
calculus used in the study of logarithmic intertwining operators and
logarithmic tensor category theory for modules for a vertex operator
algebra.  This extension has a variety of interesting arithmetic
properties.  We develop one such result here, the aforementioned
recursive identity.  We have applied this identity elsewhere to
certain formal series expansions related to a general formal Taylor
theorem and these series expansions in turn yield a sequence of
combinatorial identities which have as special cases certain classical
combinatorial identities involving (separately) the Stirling numbers
of the first and second kinds.
\end{abstract}

\maketitle

\renewcommand{\theequation}{\thesection.\arabic{equation}}
\renewcommand{\thetheorem}{\thesection.\arabic{theorem}}
\setcounter{equation}{0} \setcounter{theorem}{0}
\setcounter{section}{0}

\section{Introduction} 
In \cite{M} and \cite{HLZ} logarithmic formal calculus was used to set
up certain structure for the treatment of logarithmic intertwining
operators and ultimately logarithmic tensor category theory for
modules for a vertex operator algebra.  One particular foundational
step in \cite{HLZ} involved two expansions of certain formal series
which yielded a classical combinatorial identity involving Stirling
numbers of the first kind (see (3.17) in \cite{HLZ}), which was used
to solve a problem posed in \cite{Lu} (see Remark 3.8 in \cite{HLZ}).
These series expansions were worked out during the course of a proof
of a very general logarithmic formal Taylor theorem (see Theorem 3.6
in \cite{HLZ}).  A detailed treatment of an efficient algebraic method
to obtain formal Taylor theorems for much more general kinds of
``formal functions'' was given in \cite{R1}.  This method was
demonstrated on a space involving formal versions of {\it iterated}
logarithmic and exponential variables, extending the setting used in
\cite{HLZ}.  The method of proof bypasses any series expansions.  The
series expansions generalizing the one appearing in the proof of
Theorem 3.6 in \cite{HLZ} were carried out in \cite{R2} and they
yielded, among other identities, both the identity involving Stirling
numbers of the first kind mentioned above and also an analogous
identity involving the Stirling numbers of the second kind.  It was
during the course of that work that the recursive identity which is
the subject of this paper was discovered (and applied).  As is often
the case, the side project turned out to be quite as interesting as
the original work.  Both the statement and proof of this recursive
identity are natural and quite different in methods and philosophy
from the original work which suggested them, and we felt that they
deserved a separate treatment.  The formal calculus of iterated
logarithms and exponentials which is our setting has further
interesting arithmetic properties of which this recursive identity is
just one (e.g. see \cite{R4}).

 The purpose of this paper is to give a complete proof of this
identity, which we now state (without all the definitions) as a
preview: For $n \in \mathbb{Z}$, we have
\begin{align*}
\ell_{n+1}(x+y)=\ell_{n+1}(x)+\log
\left(1+\left(\frac{\ell_{n}(x+y)-\ell_{n}(x)}
{\ell_{n}(x)}\right)\right),
\end{align*}
where $\ell_{n}(x)$ is a formal analogue of the $(-n)$-th iterated
exponential for $n < 0$ and the n-th iterated logarithm for $n > 0$
and $\ell_{0}(x)$ is a formal analogue of $x$ itself, where for a
formal object $X$,
\begin{align*}
\log (1+X)=\sum_{i \geq 1}\frac{(-1)^{i-1}}{i}X^{i},
\end{align*}
whenever this formal expression is well defined (so that there are two
different types of ``logarithm,'' as is also the case in \cite{HLZ}).
We note also the following equivalent form of the above identity:
\begin{align*}
\ell_{n}(x+y)=\ell_{n}(x)e^{\left(\ell_{n+1}(x+y)-\ell_{n+1}(x)\right)},
\end{align*}
where for a formal object $X$
\begin{align*}
e^{X}=\sum_{i \geq 0}\frac{X^{i}}{i!},
\end{align*}
whenever this is well defined.  While the result is heuristically clear, it is
certainly not trivial when, for instance, one considers the fragile
blend of two different types of logarithm, and the
proof is hardly obvious.  This work is a slightly updated version of
part of \cite{R3}.  

This recursive identity involves expressions with a formal translation
occuring in the arguments of some of the expressions.  Analogous
identities may be developed where we use more general formal changes
of variable.  One could develop such identities using parallel
arguments to those we employ here.  However, there is a nicer way to
obtain further identities of this type using a different idea, which
is developed in \cite{R4} (see also \cite{R1}), to which we refer the
interested reader.
\section{Formal iterated logarithms and exponentials}
\label{sec:varch}
Let $\ell_{n}(x)$ be formal commuting variables for $n \in
\mathbb{Z}$. We consider the algebra with an underlying vector space
basis consisting of all elements of the form
\begin{align*}
\prod_{i \in \mathbb{Z}}\ell_{i}(x)^{r_{i}},
\end{align*}
where $r_{i} \in \mathbb{C}$ for all $i \in \mathbb{Z}$, and all but
finitely many of the exponents $r_{i}=0$.  The multiplication is the
obvious one (when multiplying two monomials simply add the
corresponding exponents and linearly extend).  We call this algebra
\begin{align*}
\mathbb{C}\{[\ell]\}.
\end{align*}
We let $\frac{d}{dx}$ be the unique derivation on $\mathbb{C}\{[\ell]\}$
satisfying
\begin{align*}
\frac{d}{dx}\ell_{-n}(x)^{r}
&=r\ell_{-n}(x)^{r-1}\prod_{i=-1}^{-n}\ell_{i}(x),\\
\frac{d}{dx}\ell_{n}(x)^{r}
&=r\ell_{n}(x)^{r-1}\prod_{i=0}^{n-1}\ell_{i}(x)^{-1},\\
\text{\rm and   } \qquad \frac{d}{dx}\ell_{0}(x)^{r}
&=r\ell_{0}(x)^{r-1},
\end{align*}
for $n > 0$ and $r \in \mathbb{C}$.  
\begin{remark} \rm
\label{rem:logexpsecret}
Secretly, $\ell_{n}(x)$ is the $(-n)$-th iterated
exponential for $n < 0$ and the $n$-th iterated logarithm for $n > 0$
and $\ell_{0}(x)$ is $x$ itself.
\end{remark}
\begin{remark} \rm
To see that this does indeed uniquely define a derivation, we note that
$\frac{d}{dx}$ must coincide with the unique linear map satisfying
\begin{align*}
\frac{d}{dx}\prod_{i \in \mathbb{Z}}\ell_{i}(x)^{r_{i}}
=\sum_{j \in \mathbb{Z}}\frac{d}{dx}\ell_{j}(x)^{r_{j}}
\prod_{i \neq j \in \mathbb{Z}}\ell_{i}(x)^{r_{i}},
\end{align*}
on a basis of $\mathbb{C}\{[\ell]\}$.  This establishes uniqueness.
We need to check that this linear map is indeed a derivation.  It is
routine and we leave it to the reader to verify that it is enough to
check that
\begin{align*}
\frac{d}{dx}(ab)=\left(\frac{d}{dx}a\right)b+a\left(\frac{d}{dx}b\right)
\end{align*}
for basis elements $a$ and $b$.  Another routine calculation reduces
the case to where $a=\ell_{i}(x)^{r}$ and $b=\ell_{i}(x)^{s}$ for $r,s
\in \mathbb{C}$.  Checking this case is trivial once one notes that
\begin{align*}
\frac{d}{dx}\ell_{j}(x)^{r}=r\ell_{j}(x)^{r-1}\frac{d}{dx}\ell_{j}(x).
\end{align*}
\end{remark}
If we let $x$ and $y$ be independent formal variables, then the formal
exponentiated derivation $e^{y\frac{d}{dx}}$, defined by the
expansion, $\sum_{k \geq 0}y^{k}\left(\frac{d}{dx}\right)^{k}/k!$,
acts on a (complex) polynomial $p(x)$ as a formal translation in $y$.
That is, as the reader may easily verify, we have
\begin{align}
e^{y\frac{d}{dx}}p(x)=p(x+y).\label{eq:PFTT}
\end{align}
This motivates the following definition (as in \cite{R1}).
\begin{defi}
Let
\begin{align*}
\ell_{n}(x+y)=e^{y\frac{d}{dx}}\ell_{n}(x) \qquad \text{for} 
\quad n \in \mathbb{Z}.
\end{align*}
\end{defi}

We establish a recursive identity for $\ell_{n}(x+y)$
in terms of $\ell_{n-1}(x+y)$ (and inversely in terms of
$\ell_{n+1}(x+y)$) for all $n \in \mathbb{Z}$.  

Our approach is based on the following identity:
\begin{align}
\label{eq:itloglim}
\text{lim}_{r \rightarrow0}\left(\left(\frac{d}{dx}\right)^{m}
\frac{\left(\ell_{n}(x)\right)^{r}}{r}\right)
=\left(\frac{d}{dx}\right)^{m}\ell_{n+1}(x)
\quad (m \geq 1).
\end{align}
But we shall need to define just what we mean by taking a limit in
this context in order for the above expression to make precise sense.
We first define a new space.
\begin{defi} \rm
We let $F(\mathbb{Z_{+}},\ell)$ be the complex vector space of functions
from the positive integers into $\mathbb{C}\{[\ell]\}$.
\end{defi} 
We may define a ``lifting'' of $\frac{d}{dx}$ on $F(\mathbb{Z_{+}},\ell)$.
\begin{defi} \rm
For $f$ and $g$ $\in F(\mathbb{Z_{+}},\ell)$, we say that
$g=\frac{d}{dx}f$ when $g(r)=\frac{d}{dx}f(r)$ for all $r \geq 0$.
\end{defi}
Of course, $\frac{d}{dx}f$ may not exist for all $f \in
F(\mathbb{Z_{+}},\ell)$.  We shall actually be interested in a
subspace of $F(\mathbb{Z_{+}},\ell)$, which we call
$P(\mathbb{Z_{+}},\ell)$ on which $\frac{d}{dx}f$ does always exist.
\begin{defi} \rm
 We let $P(\mathbb{Z_{+}},\ell)\subset F(\mathbb{Z_{+}},\ell)$ 
be the space of functions $f(r)$, from the nonzero natural numbers into
$\mathbb{C}\{[\ell]\}$
, which may be
 represented in the form
\begin{align*}
f(r)=\sum_{j \geq 0} q_{j}(r)
\prod_{i \in \mathbb{Z}} \ell_{i}(x)^{p_{i,j}(r)},
\end{align*}
where $q_{j}(r),p_{i,j}(r)$ are complex polynomials in $r$ for all $j
\geq 0$, $i \in \mathbb{Z}$ and where we further require that for all
$j \geq 0$ there exists some $I_{j} \geq 0$ such that $p_{i,j}(r)=0$
when $|i| \geq I_{j}$ and finally that there exists $J \geq 0 $ such
that $q_{j}(r)$ is the zero polynomial for $j \geq J$.  We call such a
representation a \it{formal polynomial form} \rm of the function.  The
function is given by the obvious substitution procedure for $r$ in the
formal polynomial form.
\end{defi}
\begin{defi} \rm
For $f(r) \in P(\mathbb{Z_{+}},\ell)$ we say a formal polynomial form
representation,
\begin{align*}
f(r)=
\sum_{j \geq 0} q_{j}(r)\prod_{i \in \mathbb{Z}} \ell_{i}(x)^{p_{i,j}(r)},
\end{align*}
is in \it{reduced formal polynomial form} \rm or \it{reduced form},
\rm when for all $j \neq k$, $j,k \geq 0$ there is some $i \in
\mathbb{Z}$ such that
\begin{align*}
p_{i,j}(r)  \neq  p_{i,k}(r).
\end{align*}
\end{defi}
Then we get
\begin{prop}
If $f(r)$ $\in P(\mathbb{Z_{+}},\ell)$, then it is uniquely expressible in
reduced formal polynomial form.
\end{prop}
\begin{proof}
Let us say that 
\begin{align*}
M(r)=q(r)\prod_{i \in \mathbb{Z}}\ell_{i}(x)^{p_{i}(r)}
\end{align*}
is a monomial summand in one reduced formal polynomial form of $f(r)$.
Then consider any other reduced formal polynomial form expression for
$f(r)$.  Since two formally unequal complex polynomials can only agree
for a finite number of substitution values, it is not difficult to see
that there must be a monomial summand in the second reduced polynomial
form of the form
\begin{align*}
N(r)=\bar {q} (r)\prod_{i \in \mathbb{Z}}\ell_{i}(x)^{p_{i}(r)}.
\end{align*}
But since our forms are reduced, then in fact $N(r)$ is the only
monomial summand of this form in the second representation, and
therefore $q(r)=\bar{q}(r)$.  The result now obviously follows by
induction.
\end{proof}
It is now easy to define what is meant by $\text{lim}_{r \rightarrow
0}f(r)$ when $f(r) \in P(\mathbb{Z_{+}},\ell)$.  One simply expresses
$f(r)$ in its unique reduced formal polynomial expansion and
substitutes $0$ for $r$ to yield a well-defined element of 
$\mathbb{C}\{[\ell]\}$.
Before we return to the identity which we want we should note that
$P(\mathbb{Z_{+}},\ell)$ is obviously closed under $\frac{d}{dx}$.  It
is also necessary to prove one lemma.
\begin{lemma}
For any $A_{r}(x) \in P(\mathbb{Z_{+}},\ell)$ we have that
\begin{align*}
\lim_{r \rightarrow 0}\frac{d}{dx}A_{r}(x)
=\frac{d}{dx}\lim_{r \rightarrow 0}A_{r}(x).
\end{align*}
\end{lemma}
\begin{proof}
Since $\frac{d}{dx}$ is linear we only have to consider the case where
$A_{r}(x)$ is a monomial.  For convenience we call $\lim_{r
\rightarrow 0}A_{r}(x)=A_{0}(x)$.  Let $A_{r}(x)=B_{r}(x)C_{r}(x).$
then
\begin{align*}
\lim_{r \rightarrow 0}\frac{d}{dx}A_{r}(x)&=\lim_{r \rightarrow
0}\frac{d}{dx}(B_{r}(x)C_{r}(x))\\ &=(\lim_{r \rightarrow
0}\frac{d}{dx}B_{r}(x))C_{0}(x)+B_{0}(x)\lim_{r \rightarrow 0}
\frac{d}{dx}(C_{r}(x))
\end{align*}
and
\begin{align*}
\frac{d}{dx}\lim_{r \rightarrow 0}A_{r}(x)
&=\frac{d}{dx} \lim_{r \rightarrow 0} (B_{r}(x)C_{r}(x))\\
&=(\frac{d}{dx}B_{0}(x))C_{0}(x))+B_{0}(x)\frac{d}{dx}C_{0}(x)),
\end{align*}
which means that we only need consider the case where
$A_{r}(x)=p(r)\ell_{i}(x)^{q(r)}$ where $i \in \mathbb{Z}$.  Now we
get:
\begin{align*}
\lim_{r \rightarrow 0}\frac{d}{dx}A_{r}(x)&=\lim_{r \rightarrow
0}q(r)p(r)\ell_{i}(x)^{q(r)-1}\frac{d}{dx}\ell_{i}(x)\\
&=q(0)p(0)\ell_{i}(x)^{q(0)-1}\frac{d}{dx}\ell_{i}(x)\\
&=\frac{d}{dx}p(0)\ell_{i}(x)^{q(0)}\\ 
&=\frac{d}{dx}\lim_{r \rightarrow 0}A_{r}(x).
\end{align*}
\end{proof}
We now prove the desired identity (\ref{eq:itloglim}).
\begin{lemma}
For $m \geq 1$ and $n \in \mathbb{Z}$
\begin{align*}
\lim_{r \rightarrow 0}
\left(\left(\frac{d}{dx}\right)^{m}
\frac{\left(\ell_{n}(x)\right)^{r}}{r}\right)
=\left(\frac{d}{dx}\right)^{m}\ell_{n+1}(x).
\end{align*}
\end{lemma}
\begin{proof}
First note that $\frac{d}{dx}\frac{\left(\ell_{n}(x)\right)^{r}}{r} \in
P(\mathbb{Z_{+}},\ell)$.  Then we may calculate to get:
\begin{align*}
&\underset{{r \rightarrow 0}}{\text{lim}}
\left(\left(\frac{d}{dx}\right)^{m}
\frac{\left(\ell_{n}(x)\right)^{r}}{r}\right)
=\left(\frac{d}{dx}\right)^{m-1}\text{lim}_{r \rightarrow 0}
\left(\frac{d}{dx}\frac{\left(\ell_{n}(x)\right)^{r}}{r}\right)\\
&=\left(\frac{d}{dx}\right)^{m-1}\text{lim}_{r \rightarrow 0}
\left(\ell_{n}(x)^{r-1}\frac{d}{dx}\ell_{n}(x)\right)\\
&=\left(\frac{d}{dx}\right)^{m-1}\left(\ell_{n}(x)^{-1}
\frac{d}{dx}\ell_{n}(x)\right).\\
\end{align*}
And now we proceed in the two separate cases $n \geq 0$ and $n \leq
-1$.  First, when $n \geq 0$ we have,
\begin{align*}
\underset{{r \rightarrow 0}}{\text{lim}}
\left(\left(\frac{d}{dx}\right)^{m}
\frac{\left(\ell_{n}(x)\right)^{r}}{r}\right)
&=\left(\frac{d}{dx}\right)^{m-1}
\left(\ell_{n}(x)^{-1}\prod_{i=0}^{n-1}
\left(\ell_{i}(x)\right)^{-1}\right)\\
&=\left(\frac{d}{dx}\right)^{m-1}\prod_{i=0}^{n}
\left(\ell_{i}(x)\right)^{-1}\\
&=\left(\frac{d}{dx}\right)^{m}\ell_{n+1}(x).
\end{align*}
And second, when $n \leq -1$ we have,
\begin{align*}
\underset{{r \rightarrow 0}}{\text{lim}}
\left(\left(\frac{d}{dx}\right)^{m}
\frac{\left(\ell_{n}(x)\right)^{r}}{r}\right)
&=\left(\frac{d}{dx}\right)^{m-1}
\left(\ell_{n}(x)^{-1}\prod_{i=-1}^{n}
\ell_{i}(x)\right)\\
&=\left(\frac{d}{dx}\right)^{m-1}
\left(\prod_{i=-1}^{n+1}
\ell_{i}(x)\right)\\
&=\left(\frac{d}{dx}\right)^{m}\ell_{n+1}(x).
\end{align*}
\end{proof}
With some care, we now see that for $n \in \mathbb{Z}$
\begin{align*}
&\text{lim}_{r \rightarrow 0}\left(e^{y\frac{d}{dx}}
\left(\frac{\left(\ell_{n}(x)\right)^{r}}{r}\right)-
\frac{\ell_{n}(x)^{r}}{r}\right)=e^{y\frac{d}{dx}}
\ell_{n+1}(x)-\ell_{n+1}(x).
\end{align*}
One must note that indeed
$e^{y\frac{d}{dx}}\left(\frac{\left(\ell_{n}(x)\right)^{r}}{r}\right)-
\frac{\ell_{n}(x)^{r}}{r} \in P(\mathbb{Z_{+}},\ell)$ because the
first term of $e^{y\frac{d}{dx}} \left(\frac{\left(\ell_{n}(x)
\right)^{r}}{r}\right)$ cancels $\frac{\ell_{n}(x)^{r}}{r}$.  Next
we get for $n \in \mathbb{Z}$
\begin{align*}
&\ell_{n+1}(x+y)=\ell_{n+1}(x)+ \text{lim}_{r \rightarrow
0}\left(\frac{\ell_{n}(x+y)^{r}-\ell_{n}(x)^{r}}{r}\right).
\end{align*}
But we don't want the limit in the expression, so, recalling that $r$
stands for a positive integer, we get:
\begin{align*}
&\ell_{n+1}(x+y)=\ell_{n+1}(x)+\underset{r \rightarrow
0}{\text{lim}} \left( \frac{ \left(
\ell_{n}(x)+\left(\ell_{n}(x+y)-\ell_{n}(x) \right)
\right)^{r}-\ell_{n}(x)^{r}}{r}\right)\\
&=\ell_{n+1}(x)+\underset{r \rightarrow 0}{\text{lim}} \sum_{p \geq
1}\frac{r(r-1) \cdots (r-(p-1))}{rp!}\ell_{n}(x)^{r-p} \left(
\ell_{n}(x+y)-\ell_{n}(x)\right) ^{p} \\ 
&=\ell_{n+1}(x)+\sum_{p
\geq 1}\frac{(-1)^{p-1}}{p}
\left(\frac{\ell_{n}(x+y)-\ell_{n}(x)}{\ell_{n}(x)}\right)^{p}\\
&=\ell_{n+1}(x)+\log
\left(1+\left(\frac{\ell_{n}(x+y)-\ell_{n}(x)}{\ell_{n}(x)}\right)
\right),
\end{align*}
where for a formal object $X$ 
\begin{align*}
\log (1+X)=\sum_{i \geq 1}\frac{(-1)^{i-1}}{i}X^{i},
\end{align*}
whenever this formal expression is well defined.
\begin{theorem} For $n \in \mathbb{Z}$ we have
\begin{align}
\label{eq:itlogrec}
\ell_{n+1}(x+y)=\ell_{n+1}(x)+\log
\left(1+\left(\frac{\ell_{n}(x+y)-\ell_{n}(x)}
{\ell_{n}(x)}\right)\right).
\end{align}
\end{theorem}
\begin{flushright} $\qed$ \end{flushright}

We note that we may solve for $\ell_{n}(x+y)$ in (\ref{eq:itlogrec})
to get for all $n \in \mathbb{Z}$
\begin{align}
\label{eq:itlogrec2}
\ell_{n}(x+y)=\ell_{n}(x)e^{\left(\ell_{n+1}(x+y)-\ell_{n+1}(x)\right)},
\end{align}
where for a formal object $X$
\begin{align*}
e^{X}=\sum_{i \geq 0}\frac{X^{i}}{i!},
\end{align*}
whenever this is well defined.  We used that $e^{\log (1+X)}=1+X$,
which is perhaps checked most easily by calculating that
$\frac{d}{dx}\left(e^{\log (1+X)}(1+X)^{-1}\right)=0$ which gives
$e^{\log (1+X)}=c(1+X)$ for some constant, $c$ which is, in turn,
solved for by substituting $0$ for $X$ or, in other words, checking
the constant term.
\begin{remark} \rm
Both (\ref{eq:itlogrec}) and (\ref{eq:itlogrec2}) make sense
heuristically, as may be seen easily, when one recalls that
$\ell_{n}(x)$ is ``really'' an iterated logarithm or exponential.
\end{remark}

\noindent {\small \sc Department of Mathematics, Rutgers University,
Piscataway, NJ 08854} 
\\ {\em E--mail
address}: thomasro@math.rutgers.edu
\end{document}